\documentclass[11pt,a4paper]{amsart}

\usepackage{amsmath}
\usepackage{mathtools}
\usepackage{amsthm}
\usepackage{amssymb}
\usepackage{hyperref}

\newcommand*\fullref[3][\relax]{%
  \ifdefined\hyperref%
    {\hyperref[#3]{#2\penalty 200\ \ref*{#3}#1}}%
  \else%
    {#2\penalty 200\ \relax\ref{#3}#1}%
  \fi%
}

\usepackage{tikz}
\usetikzlibrary{shapes.misc}

\numberwithin{equation}{section}

\DeclarePairedDelimiter{\parens}{\lparen}{\rparen}

\DeclarePairedDelimiter{\set}{\{}{\}}

\newcommand*\nset{\mathbb{N}}

\newtheorem*{theorem*}{Theorem}

\newcommand*{\sym}[1]{\overline{#1}}

\newcommand*{\fsa}[1]{\mathcal{#1}}

\begin{document}

\title{Every numerical semigroup arises as an automaton monoid}

\author{Tara Macalister Brough}
\address[T. Macalister Brough]{%
Center for Mathematics and Applications (NOVA Math)\\
NOVA School of Science and Technology\\
NOVA University of Lisbon\\
2829--516 Caparica\\
Portugal
}
\email{%
t.brough@fct.unl.pt
}

\thanks{For the first and second authors, this work is funded by national funds through the FCT -
  Fundação para a Ciência e a Tecnologia, I.P., under the scope of the projects UIDB/00297/2020 and UIDP/00297/2020
  (Center for Mathematics and Applications).}

\author{Alan J. Cain}
\address[A.J. Cain]{%
Center for Mathematics and Applications (NOVA Math)\\
NOVA School of Science and Technology\\
NOVA University of Lisbon\\
2829--516 Caparica\\
Portugal
}
\email{%
a.cain@fct.unl.pt
}

\author{Jan Philipp Wächter}
\address[J.P. Wächter]{%
Dipartimento di Matematica\\
Politecnico di Milano\\
Piazza Leonardo da Vinci, 32\\
20133 Milano, Italy
}

\thanks{At the time the research leading to this paper was initiated, the third author was affiliated to Centro de
  Matemática da Universidade do Porto (CMUP), which is financed by national funds through FCT – Fundação para a Ciência
  e Tecnologia, I.P., under the project with reference UIDB/00144/2020. Currently, he is funded by the Deutsche
  Forschungsgemeinschaft (DFG, German Research Foundation) – 492814705.}

\thanks{This work was partially supported by the Fundação para a Ciência e a Tecnologia (Portuguese Foundation for
  Science and Technology) through the project PTDC/MAT-PUR/31174/2017.}

\begin{abstract}
  This paper shows how to construct explicitly an automaton that generates an arbitrary numerical semigroup.
\end{abstract}


\maketitle

\section{Introduction}

The aim of this paper is to give a constructive proof of the theorem stated in the title:

\begin{theorem*}
  Every numerical semigroup arises as an automaton monoid.
\end{theorem*}

Let $\nset_0$ denote the set of natural numbers including zero. A numerical semigroup is a submonoid of
$\parens{\nset_0,+}$ with finite complement; note that a numerical semigroup necessarily contains $0$. The theory of
numerical semigroups is ultimately rooted in the Frobenius problem of finding the largest natural number not expressible
as a sum of given natural numbers. In the twentieth century, numerical semigroups gained in importance because of their
applications in algebraic geometry and connections to number theory. For background on numerical semigroups, see, for
example, \cite{rosales_numerical}.

Automaton semigroups are semigroups generated by the actions of Mealy machines (that is, finite automata with output).
They emerged as a generalization of automaton groups, which arose from the construction of `exotic' groups such as
infinite torsion groups (for background, see, for example, \cite{nekrashevych_self}). A substantial theory of automaton
semigroups has emerged over the last fifteen years (see, for example,
\cite{maltcev_cayley,bc_2auto,dangeli_onthestructure}). An automaton monoid is simply a monoid generated as a semigroup
by an automaton.

A particular topic of research within this theory has concerned showing that certain (semi)groups arise or do not arise
as automaton (semi)groups. The theorem above contributes to this topic by showing that all semigroups in a very natural
and important class are automaton semigroups. The proof in \fullref{Section}{sec:proof} can be used to construct as an
automaton semigroup any numerical semigroup given by its finite complement in $\nset_0$.

The theorem is also important because it extends a contrast between $(\nset,+)$ (which is \emph{not} an automaton
semigroup \cite[Proposition~4.3]{c_1auto}) and $(\nset_0,+)$ (which \emph{is}). The first and second authors proved
previously that no subsemigroup of $(\nset,+)$ is an automaton semigroup; indeed, no non-trivial subsemigroup of this
semigroup with an adjoined multiplicative zero is an automaton semigroup \cite[Theorem~15]{bc_2auto}. (Note that the
adjoined zero here is an \emph{absorbing} element $z$ with $xz = zx =z$, \emph{not} an adjoined additive zero or neutral
element.) As a consequence, although each numerical semigroup arises as an automaton semigroup, the semigroup formed by
removing the identity (that is, the neutral element $0 \in \nset_0$) from a numerical semigroup is never an automaton
semiroup.

\section{Notation \& background references}

%This section sets up notation and gives references for the necessary concepts from the theory of automaton semigroups.
%TB: I feel the section is so short that this sentence is rather unnecessary.

For any symbol $a$, denote by $a^\omega$ the infinite sequence of symbols $a$.
For a set $A$, denote by $A^\omega$ the set of all infinite sequences of symbols from $A$.

For the definition of automata (Mealy machines) and their actions on sequences from some symbol alphabet, see
\cite[\S~2]{c_1auto}. The semigroup generated by an automaton $\fsa{M}$ is denoted $\Sigma\parens{\fsa{M}}$.

We use the convention that states of an automaton act from the right.  If $w$ is a word in the states of an automaton,
and $\alpha$ an input string, then we write $\alpha\cdot w$ for the action of $w$ on $\alpha$.

For the definition of wreath recursions and how to calculate using them, see \cite[\S~3.1]{c_1auto}.

\section{Proof of the theorem}
\label{sec:proof}

Let $S$ be a numerical semigroup. Then $S = \nset_0 \setminus F$ for some finite $F \subseteq \nset$. The strategy of
the proof is first to build a special automaton (different from the standard automaton for $\nset_0$) that generates
$\nset_0$ as an automaton semigroup, and where $0$ has an identical action. In this automaton, the state set includes
$0$ and all elements of an initial interval of the natural numbers that includes $F$. The automaton is such that the
states in $F$ can be deleted to leave an automaton that generates $S$.

Let $n \in \nset$. Let $\fsa{M}_n$ be the automaton with state set $\set{0,1,n}$ operating on the symbol alphabet
$A = \set{\sym{0},\sym{1},\ldots,\sym{n}}$, with transitions as defined as follows:
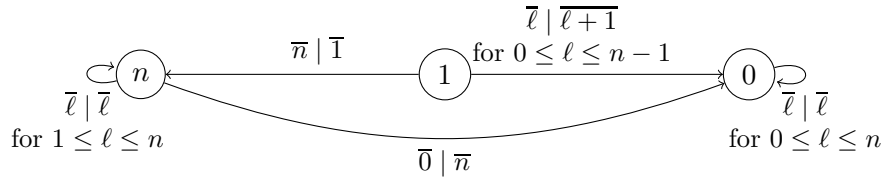
\[
  \begin{tikzpicture}
    \begin{scope}[
      every node/.style={
        draw,
        circle,
      }
      ]
      \node (n) at (0,0) {$n$};
      \node (1) at (4,0) {$1$};
      \node (0) at (8,0) {$0$};
    \end{scope}
    \begin{scope}[
      every node/.style={
        node font=\small,
        align=center,
        }
      ]

      \draw[->,loop right] (0) edge node[below,yshift=-1mm] {$\sym{\ell}\mid\sym{\ell}$\\for $0 \leq \ell \leq n$} (0);
      \draw[->] (1) edge node[above,pos=.4] {$\sym{\ell} \mid \sym{\ell+1}$\\for $0 \leq \ell \leq n-1$} (0);
      \draw[->] (1) edge node[above,pos=.4] {$\sym{n} \mid \sym{1}$} (n);
      \draw[->,loop left] (n) edge node[below,yshift=-1mm] {$\sym{\ell}\mid\sym{\ell}$\\for $1 \leq \ell \leq n$} (n);
      \draw[->] (n) edge[bend right=20] node[below] {$\sym{0} \mid \sym{n}$} (0);
    \end{scope}
  \end{tikzpicture}
\]
In terms of wreath recursions, the elements $0$, $1$, and $n$ are given by
\begin{align*}
  0       & = \parens{0,0,0,\ldots,0,0}
      \begin{pmatrix}
        \sym{0} & \sym{1} & \sym{2} & \cdots & \sym{n-1} & \sym{n} \\
        \sym{0} & \sym{1} & \sym{2} & \cdots & \sym{n-1} & \sym{n} \\
      \end{pmatrix};                  \displaybreak[0]\\
  1       & = \parens{0,0,0,\ldots,0,n}
      \begin{pmatrix}
        \sym{0} & \sym{1} & \sym{2} & \cdots & \sym{n-1} & \sym{n} \\
        \sym{1} & \sym{2} & \sym{3} & \cdots & \sym{n}   & \sym{1} \\
      \end{pmatrix};                  \displaybreak[0]\\
  n       & = \parens{0,n,n,\ldots,n,n}
      \begin{pmatrix}
        \sym{0} & \sym{1} & \sym{2} & \cdots & \sym{n-1} & \sym{n} \\
        \sym{n} & \sym{1} & \sym{2} & \cdots & \sym{n-1} & \sym{n} \\
      \end{pmatrix}.
\end{align*}

Note that $0$ acts identically on all sequences in $A^\omega$ and so is an identity for $\Sigma\parens{\fsa{M}_n}$,
which is therefore a monoid. The semigoup $\Sigma\parens{\fsa{M}_n}$ will turn out to be $(\nset_0,+)$ and the elements
$0$, $1$, and $n$ have their natural interpretation, but for the purposes of this proof, products are written
multiplicatively, not additively.

Calculation using wreath recursions, and using the fact that $0$ is an identity for $\Sigma\parens{\fsa{M}_n}$, shows
that for any $m \in \set{1,\ldots,n}$,
\begin{equation}
  \label{eq:powers-of-1}
  \begin{aligned}
    1^m  ={} & \parens{\overbrace{\strut 0,0,\ldots,0}^{\mathclap{\substack{n-m+1\\\text{components}}}},\overbrace{\strut n,\ldots n,n}^{\mathclap{\substack{m\\\text{components}}}}} \\
             & \;\;
        \begin{pmatrix}
          \sym{0}  & \sym{1}   & \cdots & \sym{n-m} & \sym{n-m+1} & \cdots & \sym{n-1} & \sym{n} \\
          \sym{m}  & \sym{m+1} & \cdots & \sym{n}   & \sym{1}     & \cdots & \sym{m-1} & \sym{m} \\
        \end{pmatrix}.
  \end{aligned}
\end{equation}
In particular, \eqref{eq:powers-of-1} shows that
\[
  1^n ={} \parens{0,\overbrace{\strut n,\ldots n,n}^{\mathclap{\substack{n\\\text{components}}}}}
        \begin{pmatrix}
          \sym{0}  & \sym{1}   & \cdots & \sym{n-1} & \sym{n} \\
          \sym{n}  & \sym{1} & \cdots  & \sym{n-1} & \sym{n}  \\
        \end{pmatrix} = n.
\]
Hence $\Sigma\parens{\fsa{M}_n}$ is generated as a monoid by the element $1$.

The next step is to show that powers of $1$ are all distinct in $\Sigma\parens{\fsa{M}_n}$. To prove this, let
$m \in \nset$ and consider the action of $1^m$ on the sequence $\sym{0}^\omega$. Suppose that $m = qn + r$, where
$q,r \in \nset$ with $0 \leq r \leq n-1$.
\begin{align*}
  &\sym{0}^\omega\cdot 1^m \\
  ={}&\sym{1}\,\sym{0}^\omega \cdot 1^{m-1}\displaybreak[0]\\
  ={}&\sym{2}\,\sym{0}^\omega \cdot 1^{m-2}\displaybreak[0]\\
  & \vdots \displaybreak[0]\\
  ={}&\sym{n}\,\sym{0}^\omega\cdot 1^{m-n} \displaybreak[0]\\
  ={}&\sym{1}\parens[\big]{\sym{0}^\omega \cdot n} \cdot 1^{m-n-1} = \sym{1}\,\sym{n}\,\sym{0}^\omega \cdot 1^{m-n-1}\\
  & \vdots \displaybreak[0]\\
  ={}&\sym{n}^2\,\sym{0}^\omega \cdot 1^{m-2n} \displaybreak[0]\\
  ={}&\sym{1}\parens[\big]{\sym{n}\,\sym{0}^\omega \cdot n} \cdot 1^{m-2n-1} = \sym{1}\,\sym{n}^2\,\sym{0}^\omega \cdot 1^{m-2n-1} \displaybreak[0]\\
  & \vdots \displaybreak[0]\\
  ={}&\sym{n}^q\,\sym{0}^\omega  \cdot 1^{m-qn} \displaybreak[0]\\
  ={}&\sym{1}\parens[\big]{\sym{n}^{q-1}\,\sym{0}^\omega \cdot n} \cdot  1^{m-qn-1}  = \sym{1}\,\sym{n}^q\,\sym{0}^\omega \cdot 1^{m-qn-1}  \\
  & \vdots \\
  ={}& \sym{r}\,\sym{n}^q\,\sym{0}^\omega.
\end{align*}
The numbers $q$ and $r$ can be read off from the result of acting on $\sym{0}^\omega$ by $1^m$. Hence the exponent
$m = qn + r$ can be recovered from the action of $1^m$. Thus $1^m = 1^\ell$ in $\Sigma\parens{\fsa{M}_n}$ if and only if
$m = \ell$

Hence $\Sigma\parens{\fsa{M}_n}$ is the free monogenic monoid generated by $1$. That is, $\Sigma\parens{\fsa{M}_n}$ is
isomorphic to $\parens{\nset_0,+}$.

Using \eqref{eq:powers-of-1}, add all elements $m \in \set{2,\ldots,n-1}$ to $\fsa{M}_n$ as redundant generators to
obtain a new automaton $\fsa{N}_n$ with $\Sigma\parens{\fsa{N}_n} = \Sigma\parens{\fsa{M}_n}$:
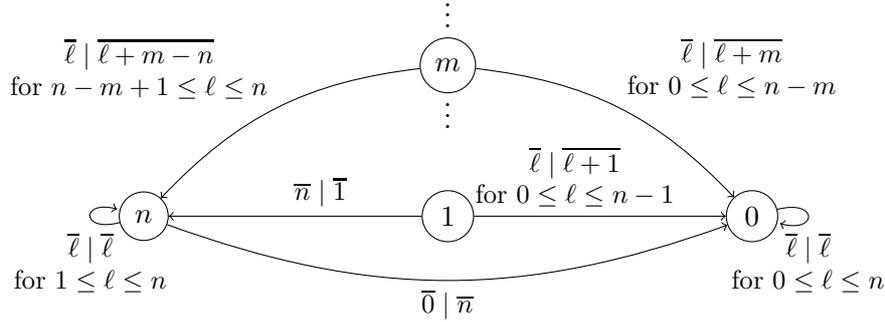
\[
  \begin{tikzpicture}
    \begin{scope}[
      every node/.style={
        draw,
        circle,
      }
      ]
      \node (n) at (0,0) {$n$};
      \node (1) at (4,0) {$1$};
      \node (0) at (8,0) {$0$};

      \node (m) at (4,2) {$m$};
    \end{scope}
    \node[anchor=south,yscale=-1] at (m.south) {$\vdots$};
    \node[anchor=south] at (m.north) {$\vdots$};
    \begin{scope}[
      every node/.style={
        node font=\small,
        align=center,
        }
      ]
      \draw[->,loop right] (0) edge node[below,yshift=-1mm] {$\sym{\ell}\mid\sym{\ell}$\\for $0 \leq \ell \leq n $} (0);
      \draw[->] (1) edge node[above,pos=.4] {$\sym{\ell} \mid \sym{\ell+1}$\\for $0 \leq \ell \leq n-1$} (0);
      \draw[->] (1) edge node[above,pos=.4] {$\sym{n} \mid \sym{1}$} (n);
      \draw[->] (m) edge[bend left=20] node[above right] {$\sym{\ell} \mid \sym{\ell+m}$\\for $0 \leq \ell \leq n-m$} (0);
      \draw[->] (m) edge[bend right=20] node[above left] {$\sym{\ell} \mid \sym{\ell+m-n}$\\for $n-m+1 \leq \ell \leq n$} (n);
      \draw[->,loop left] (n) edge node[below,yshift=-1mm] {$\sym{\ell}\mid\sym{\ell}$\\for $1 \leq \ell \leq n$} (n);
      \draw[->] (n) edge[bend right=20] node[below] {$\sym{0} \mid \sym{n}$} (0);
    \end{scope}
  \end{tikzpicture}
\]

Let $Q = \set{0,1,\ldots,n}$ and $R = \set{1,\ldots,n-1} \subseteq Q$. Note that $Q$ is the state set of $\fsa{N}_n$ and
$R$ comprises the states that have no incoming edges. Let $P \subseteq R$. Deleting every state $m \in P$ from
$\fsa{N}_n$ yields an automaton $\fsa{N}_n^{\parens{P}}$ whose action is a restriction of the action of $\fsa{N}_n$.
Hence $\Sigma\parens{\fsa{N}_n^{\parens{P}}}$ is the submonoid of
$\Sigma\parens[\big]{\fsa{N}_n} \simeq \parens{\nset_0,+}$ generated by $Q \setminus P$.

Recall that $S$ is a numerical semigroup with $S = \nset_0 \setminus F$ for some finite $F \subseteq \nset$. Let
$k \in \nset$ be such that $\set{k,k+1,k+2,\ldots} \subseteq S$; thus $F \subseteq \set{1,\ldots,k-1}$. Then
$\Sigma\parens[\big]{\fsa{N}_{2k}^{\parens{F}}}$ comprises the submonoid of $\nset_0$ generated by
$\set{1,\ldots,2k} \setminus F$. Since $F \subseteq \set{1,\ldots,k-1}$ and every element of $\nset_0$ greater than $2k$
is generated by $\set{k,\ldots,2k}$, it follows that
\begin{align*}
  \Sigma\parens[\big]{\fsa{N}_{2k}^{\parens{F}}} & = \set{0} \cup \parens[big]{\set{1,\ldots,k-1} \setminus F} \cup \set{k,\ldots,2k} \cup \set{2k+1,2k+2,\ldots} \\
                                                 & = \nset_0 \setminus F = S.
\end{align*}

\bibliography{\jobname}
\bibliographystyle{alphaabbrv}

\end{document}